# *The mixed virtual element method for grids with curved interfaces*

Franco Dassi*, Alessio Fumagalli**, Davide Losapio**, Stefano Scialò***, Anna Scotti**, Giuseppe Vacca*

\* Università degli Studi di Milano Bicocca - Member of INdAM-GNCS research group
\*\* Politecnico di Milano - Member of INdAM-GNCS research group
\*\*\* Politecnico di Torino - Member of INdAM-GNCS research group

## Abstract

In many applications the accurate representation of the computational domain is a key factor to obtain reliable and effective numerical solutions. Curved interfaces, which might be internal, related to physical data, or portions of the physical boundary, are often met in real applications. However, they are often approximated leading to a geometrical error that might become dominant and deteriorate the quality of the results. Underground problems often involve the motion of fluids where the fundamental governing equation is the Darcy law. High quality velocity fields are of paramount importance for the successful subsequent coupling with other physical phenomena such as transport. The virtual element method, as solution scheme, is known to be applicable in problems whose discretizations requires cells of general shape, and the mixed formulation is here preferred to obtain accurate velocity fields. To overcome the issues associated to the complex geometries and, at the same time, retaining the quality of the solutions, we present here the virtual element method to solve the Darcy problem, in mixed form, in presence of curved interfaces in two and three dimensions. The numerical scheme is presented in detail explaining the discrete setting with a focus on the treatment of curved interfaces. Examples, inspired from industrial applications, are presented showing the validity of the proposed approach.

## Introduction

The simulation of underground flows is a very active research field, mainly driven by the demand from the industry of increasingly robust and accurate simulation tools. The subsoil is, indeed, a porous medium characterized by complex geometries, strong heterogeneities of the properties and with an intrinsic multi-scale nature, for the presence of geometrical features and interfaces with high aspect ratios (such as fractures or channels) and forming intersections spanning several orders of magnitude. This makes underground flow simulations extremely challenging. In addition to this, simulation data are usually only available as probability distribution functions, such that many simulations are typically required to derive statistics on the desired quantities of interest. Uncertainty quantification strategies [ (Lu 2016), (Pieraccini 2020)] can be adopted to reduce the number of simulations, but the success of such tools depends, in turn, on the efficiency, robustness, reliability and accuracy of the underlying simulation method. Indeed, simulation settings derived from a realization of a random set of parameters can display complex geometries or sharp variations of parameters and the numerical tool is expected to provide an accurate result at a reasonable computational cost without requiring ad-hoc adjustments of the input data, which would alter the probability distribution of the inputs.

Discretization schemes based on the use of polygonal/polyhedral meshes are gaining an increasing popularity for underground flow simulations, especially for the possibility of an easy meshing of complex geometrical features, such as interfaces, inclusions, domain boundaries. Next to the most well-established schemes, as multi-point flux approximation [ (Sandve 2012), (Faille 2016)], new methods have been recently applied to simulations in porous media. The Virtual Element Method (VEM) [ (L. Beirao da Veiga, F. Brezzi, et al. 2015)], taking inspiration from the Mimetic Finite Differences approach [ (Lipnikov 2014), (O. S. Al-Hinai 2015), (O. W. Al-Hinai 2017), (Antonietti, et al. 2016)], allows the use of meshes of almost arbitrarily shaped elements in a finite element-like setting. The scheme is based on the definition of a discrete approximation space of virtual functions, that are never explicitly computed; suitable projectors, from the space of virtual functions to the space of polynomials, are then built, allowing to compute the discrete solution. The computability of the projectors and of the resulting discrete operators by means of the chosen set of degrees of freedom is one of the key aspects of the method. The usability of the method in the context of the simulation of underground phenomena appeared evident since its earliest appearance, for the description of the flow in poro-fractured domains [ (Benedetto, et al. 2014), (Andersen, Nilsen and Raynaud 2017), (Fumagalli and Keilegavlen 2019) (A. Fumagalli 2018), (Coulet, et al. 2020), (Mazzia, et al. 2020), (Borio, Fumagalli and Scialò 2020)]. The robustness of the VEM to highly distorted and elongated elements [ (Mascotto 2018), (Berrone and and Borio 2017)] allows for an easy meshing process also in presence of multiple intersecting interfaces, without requiring modifications of the geometry of the domain. Next to VEM and MFD, other polygonal/polyhedral methods vave been proposed, such as gradient schemes [Brenner2016], Discontinuous Galerkin approaches [ (Antonietti, et al. 2019)] and hybrid high order (HHO) methods [ (Chave, Di Pietro and and Formaggia 2018)].

The present work deals with the accurate description of single-phase flow problems in mixed form in underground reservoirs via the virtual element method in mixed formulation (MVEM) [ (L. Beirao da Veiga, et al. 2016), (Benedetto, Borio and Scialò 2017), (Dassi and Vacca 2020), (Dassi and Scacchi 2020)]. The focus is on high order approximation of complex geometries in two and three dimensions. To this aim, an extension of the MVEM to include the handling of polygonal/polyhedral elements with curvilinear boundaries [ (Beirao da Veiga, Russo and Vacca 2019), (Bertoluzza, Pennacchio and Prada 2019), (Dassi, Fumagalli, et al. 2020)] is proposed and analysed. The approach yields a discretization scheme capable of high accuracy also in presence of curvilinear features in the simulation domain. Indeed, a curved geometry can be exactly reproduced by the method, thus providing optimal error convergence trends also for high polynomial accuracy levels, since the error is not bounded by the description of the geometry, unlike classical approximation schemes where the curvilinear features of the domains are approximated by rectilinear/planar objects. Further, the use of a mixed formulation, where both pressures and fluxes are directly computed, gives local mass conservation properties to the method. The availability of a high order approximation scheme can help in the accurate description of the solution in near-well regions of the domain, or for simulations with curvilinear interfaces or domain boundaries, or for corner-point grids. Here the extension of the MVEM for curvilinear edge elements in three dimensional problems is described, and 2D/3D applications of the method to test cases of practical interest are proposed and discussed to show the applicability of this approach to underground flow simulations.

The manuscript is organized as follows: the mathematical model is briefly discussed in the first section, followed by the presentation of the discrete problem. Numerical results are shown in a third section and finally some conclusions are proposed.

## Mathematical model

In order to present the model at the basis of the proposed approach, let us consider a porous medium filled with a single incompressible fluid phase, occupying a possibly curved domain $\Omega \subset \mathbb{R}^n$, for $n = 2$ or 3, with Lipschitz continuous boundary and external unit normal $\boldsymbol{n}$. The boundary of $\Omega$, denoted as $\partial \Omega$, is divided into a portion $\partial_n \Omega$, where natural (pressure) boundary conditions are enforced, and a part $\partial_e \Omega$ with essential (flux) boundary conditions. We have the decomposition $\overline{\partial \Omega} = \overline{\partial_e \Omega} \cup \overline{\partial_n \Omega}$ and $\partial°_e \Omega \cap \partial°_n \Omega = \emptyset$, being $\partial°_n \Omega \neq \emptyset$ for solvability purposes. We further assume that $\partial_e \Omega$ and $\partial_n \Omega$ can be each divided into a finite number of regular curves:

$$\partial_n \Omega = \cup_{i=1,\dots,N_n} \Gamma_i, \qquad \partial_e \Omega = \cup_{i=N_n+1,\dots,N} \Gamma_i$$

with $\Gamma_i$ of class $C^{m+1}$, $m \geq 0$, such that, for each $i = 1, \dots, N$ it exists an invertible $C^{m+1}$ parametrization $\gamma_i: I_i \to \Gamma_i$, being $I_i \subset \mathbb{R}$ a closed interval. We will detail in the sequel the assumption on $\gamma_i$. This assumption is required for the treatment of the curvilinear boundary, following the approach proposed by (Beirao da Veiga, Russo and Vacca 2019). The same assumption also applies to curvilinear internal interfaces, which are treated in the same manner.

The porous medium is composed by heterogeneous rocks or sedimentary materials whose permeability tensor is indicated by $\boldsymbol{K}$. We assume that $\boldsymbol{K}$ is symmetric and positive definite but may change abruptly by several order of magnitudes from region to region. The fluid phase is characterized by a dynamic viscosity $\mu$, which can be considered for our purpose as a positive real number. A scalar source or sink term $f$ is considered to model fluid injection/production, e.g. through a well.

Boundary conditions are assigned on $\partial \Omega$, prescribing the natural boundary condition $\bar{p}$ and the essential boundary condition $\bar{q}$, both as given data, and thus the equations describing the Darcy velocity $q$ and fluid pressure $p$ are: *Find $(\boldsymbol{q}, p)$ such that*

$$\begin{cases} \mu \boldsymbol{q} + \boldsymbol{K} \nabla p = \boldsymbol{0} \\ \nabla \cdot \boldsymbol{q} + f = 0 \end{cases} in \ \Omega,$$

*with boundary conditions*

$$\begin{cases} p = \bar{p}, & on \ \partial_n \Omega, \\ \boldsymbol{q} \cdot \boldsymbol{n} = \bar{q}, & on \ \partial_e \Omega. \end{cases}$$

In these equations, symbols "$\nabla$" and "$\nabla \cdot$" denote the gradient and the divergence operator, respectively.

By introducing appropriate functional spaces, the previous problem can be re-written in weak formulation, suitable for its numerical resolution with the MVEM approach. In what follows we assume, for simplicity, that $\bar{q} = 0$ otherwise a lifting technique should be adopted. The first bilinear form associated to the problem is the following:

$$a(\cdot,\cdot): \boldsymbol{V}(\Omega) \times \boldsymbol{V}(\Omega) \to \mathbb{R}, \quad a(\boldsymbol{u}, \boldsymbol{v}) \coloneqq (\mu \boldsymbol{K}^{-1} \boldsymbol{u}, \boldsymbol{v})_\Omega, \quad \forall (\boldsymbol{u}, \boldsymbol{v}) \in \boldsymbol{V}(\Omega) \times \boldsymbol{V}(\Omega),$$

being $\boldsymbol{V}(\Omega)$ the space of vector functions defined as

$$\boldsymbol{V}(\Omega) \coloneqq \{\boldsymbol{v} \in [\mathrm{L}^2(\Omega)]^n : \nabla \cdot \boldsymbol{v} \in \mathrm{L}^2(\Omega), \boldsymbol{v} \cdot \boldsymbol{n} = 0 \text{ on } \partial_e \Omega\},$$

and $(\cdot,\cdot)_\omega$ the standard $\mathrm{L}^2(\omega)$-inner product. The second bilinear form is the following:

$$b(\cdot,\cdot): \boldsymbol{V}(\Omega) \times Q(\Omega) \to \mathbb{R}, \quad b(\boldsymbol{u}, v) \coloneqq -(\nabla \cdot \boldsymbol{u}, v)_\Omega, \quad \forall (\boldsymbol{u}, v) \in \boldsymbol{V}(\Omega) \times Q(\Omega),$$

with $Q(\Omega) \coloneqq \mathrm{L}^2(\Omega)$. The weak formulation of the Darcy problem then reads: *Find $(\boldsymbol{q}, p) \in \boldsymbol{V}(\Omega) \times Q(\Omega)$ such that*

$$\begin{cases} a(\boldsymbol{q}, \boldsymbol{v}) + b(\boldsymbol{v}, p) = -(\bar{p}, \boldsymbol{v} \cdot \boldsymbol{n})_{\partial_n \Omega}, & \forall \boldsymbol{v} \in \boldsymbol{V}(\Omega) \\ b(\boldsymbol{q}, v) = (f, v)_\Omega, & \forall v \in Q(\Omega) \end{cases} \tag{1}$$

Following (D. Boffi 2013) (Raviart and Thomas 1977) it is possible to show that problem (1) admits a unique solution which continuously depends from the given data. Formulation (1) of the Darcy problem is valid for all the dimensions $n$ of the ambient space $\mathbb{R}^n$.

## Discretization

In this section we introduce the numerical approximation of model (1). Computability is a key concept in the VEM technology: contrary to the classical finite element method, we do not prescribe a-priori the shape of the basis function for each element. This flexibility comes with a cost, which makes some of the local discretization matrices not directly computable. We need to introduce additional tools, in this case a suitable projection operator, to perform such computations.

As mentioned before, the boundary of $\Omega$ might be curved or, due to the problem data, internal (curved) interfaces might be present also inside the domain. We define $\Omega_h$ to be the computational grid made of $N_E$ mutually disjoint elements $E$ such that $\bigcup_{i=1,\ldots,N_E} E_i = \Omega_h$ and we require that the boundary of $\Omega$ is perfectly represented as well as possible internal interfaces. We thus have a collection of elements that might have one or more curved edges (in 2d) or faces (in 3d), with the requirement that they can be, at most, a union of a finite number of star-shaped sub-elements (L. Beirao da Veiga, et al. 2016). Grids normally used in the industry, like corner-point grids, fulfill this requirement. We name $e$ a generic edge (in 2d) or face (in 3d) and the collection of all of them as $\mathcal{E}_h$. For a single grid element $E$ its edges or faces are named $\mathcal{E}_h(E)$.

To discretize problem (1), we follow this step by step approach: *(i)* in Subsection "Discrete spaces", we construct the finite dimensional spaces $V_k(\Omega_h)$ and $Q_k(\Omega_h)$, approximations of $V(\Omega)$ and $Q(\Omega)$, respectively, both of accuracy degree $k \geq 0$. *(ii)* In Subsection "Degrees of freedom", we introduce the set of degrees of freedom to make the objects involved computable, and *(iii)* in Subsection "Projection operators and approximation", we introduce the projection operator and discuss the local approximations. In the discussion, we accurately present how to handle curved edges or faces for the grid elements.

As common for VEM, we concentrate this discussion on a single element $E \in \Omega_h$. The global approximation spaces are also introduced as "a collection" of the local ones. The latter concept will be explained in the sequel. For an element $E$, we name the local spaces as $V_k(E)$ and $Q_k(E)$. For any curved edge or face $e$ of $E$, we assume that we know the exact parametrization of it, so that we can introduce the map $\gamma_e: \mathfrak{e} \to e$, with $\mathfrak{e}$ being a straight reference interval defined on the abscissa if the ambient space is $n = 2$ or a reference planar polygon in the $uv$ plane if $n = 3$. Furthermore, we also assume that it exists the inverse map $\gamma_e^{-1}: e \to \mathfrak{e}$, deemed to be regular enough. See Figure 1 for a graphical representation.

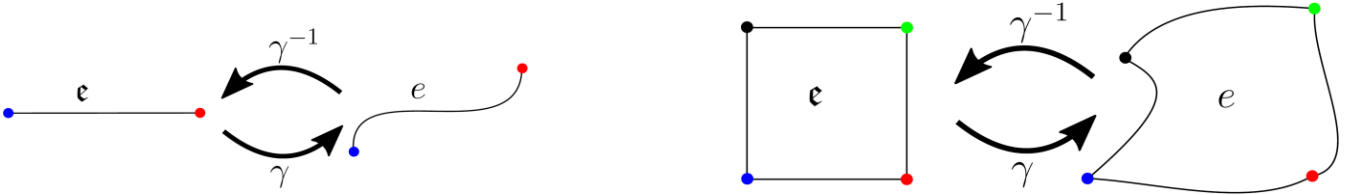

*Figure 1: Graphical example of mapping between the reference straight $\mathfrak{e}$ and the mapped curved $e$ element.*

### Discrete spaces

We first introduce some general spaces which will be used in the actual definition of $V_k(E)$ and $Q_k(E)$. We indicate with $\mathbb{P}_d(A)$ the set of polynomial on the set $A \subset \mathbb{R}^n$ of degree $\leq d$, and with $[\mathbb{P}_d(A)]^n$ the set of vector valued polynomials where each component has degree $\leq d$. We indicate with $\pi_d^n$ the dimension of $\mathbb{P}_d(A)$, which expresses the number of degrees of freedom needed to fully define a polynomial. In general we have:

$$\pi_d^n := \dim \mathbb{P}_d(A) = \frac{1}{n!} \prod_{i=1}^{n} (d + i).$$

Since a grid element may be of arbitrary shape, we cannot use a Lagrangian construction to build the polynomials. For this reason, we make use of the following scaled monomial spaces: given a straight segment or planar polygon $\mathfrak{e}$ such that $\gamma_e: \mathfrak{e} \to e$, like in Figure 1, we define

$$\mathcal{M}_d(\mathfrak{e}) := \left\{ \left( \frac{x_i - x_i(\mathfrak{e})}{h(\mathfrak{e})} \right)^j \text{ for } i \leq n \text{ and } j \leq d \right\} \quad \text{and} \quad \widetilde{\mathcal{M}}_d(e) := \{\widetilde{m} = m \circ \gamma_e^{-1} : m \in \mathcal{M}_d(\mathfrak{e})\},$$

where $x_i(\mathfrak{e})$ and $h(\mathfrak{e})$ are the $i$-th coordinate of the centre and size of $\mathfrak{e}$, respectively. The space $\mathcal{M}_d(\mathfrak{e})$ forms a base for $\mathbb{P}_d(\mathfrak{e})$, thus all the elements of the latter can be expressed as linear combination of elements in the former. Because of this, we write that $\mathbb{P}_d(\mathfrak{e}) = span(\mathcal{M}_d(\mathfrak{e}))$. Let us note that the space $\widetilde{\mathcal{M}}_d$ is not in general a polynomial space, but it will be used in the construction of the discrete problem, in particular for curved edges or faces. In the space $\widetilde{\mathcal{M}}_d$, the monomials of $\mathcal{M}_d(\mathfrak{e})$ are mapped to the actual (physical) edge or face $e$ via the map $\gamma_e^{-1}$. In the same way, we set the mapped polynomial space as

$$\widetilde{\mathbb{P}}_d(e) := \{\widetilde{p} = p \circ \gamma_e^{-1} : p \in \mathbb{P}_d(\mathfrak{e})\} = span\left( \widetilde{\mathcal{M}}_d(e) \right).$$

Similarly, we introduce for each element the scaled scalar and vector monomial spaces, respectively, as

$$\mathcal{M}_d(E) := \left\{ \left( \frac{x_i - x_i(E)}{h(E)} \right)^j \text{ for } i \leq n \text{ and } |j| \leq d \right\},$$

$$[\mathcal{M}_d(E)]^n := \{m\boldsymbol{e}_i \text{ for } i \leq n \text{ and } m \in \mathcal{M}_d(E)\},$$

with $x_i(E)$ and $h(E)$ are the $i$-th coordinate of the centre and size of $E$, respectively, $\boldsymbol{j}$ is a multi-index, and $\boldsymbol{e}_i$ is the canonical basis of $\mathbb{R}^n$. In the case the centre or the size of $E$ are not well defined or complex to compute, at discrete level it is possible to consider a "representative" centre and size that make these monomials scale in the correct way. Also in this case, the space $\mathcal{M}_d(E)$ is a basis for $\mathbb{P}_d(E)$, meaning that $\mathbb{P}_d(E) = span(\mathcal{M}_d(E))$, and similarly for $[\mathbb{P}_d(E)]^n = span([\mathcal{M}_d(E)]^n)$. In the sequel, for the numerical approximation, we will consider the monomials since they can be defined regardless of the actual shape of the element $E$ or edge/face $e$. It is possible, both for $n = 2$ or 3, to divide the space $[\mathbb{P}_d(E)]^n$ as a sum of the gradient of a polynomial in $\mathbb{P}_{d+1}(E)$ and a "remaining part". This remark will be useful in the definition of the degrees of freedom for the local space $\boldsymbol{V}_k(E)$. For $n = 2$ it is easy to define the aforementioned "remaining part", which is a vector valued polynomial indicated as $\boldsymbol{m}_d^\oplus(E)$ and given by the "rotated" linear scaled vector monomial $\boldsymbol{m}_1^\perp(E)$ multiplied by a monomial of degree $d - 1$. It is given by

$$\boldsymbol{m}_d^\oplus(E) = \boldsymbol{m}_1^\perp m_{d-1} \text{ with } \boldsymbol{m}_1^\perp(x,y) := \left(\frac{y - y(E)}{h(E)}, -\frac{x - x(E)}{h(E)}\right),$$

In the case of $n = 3$, it is possible to explicitly write it, but it will become trickier and it is out of the scope of the present work. For this, we refer to the article (Dassi and Scacchi 2020). We indicate the space of all $\boldsymbol{m}_d^\oplus(E)$ as $\mathcal{G}_d^\oplus(E)$ so that the following relation can be written:

$$[\mathbb{P}_d(E)]^n = \{\nabla p : p \in \mathbb{P}_{d+1}(E)\} \cup span\left(\mathcal{G}_d^\oplus(E)\right).$$

We introduce the local discrete space for the vector trial $\boldsymbol{q}$ and test $\boldsymbol{v}$ fields of problem (1). We have the following approximation space

$$\boldsymbol{V}_k(E) := \{\boldsymbol{v} \in V(E) : \boldsymbol{v} \cdot \boldsymbol{n} \in \widetilde{\mathbb{P}}_k(e) \; \forall e \in \mathcal{E}_h(E) \text{ and } \nabla \cdot \boldsymbol{v} \in \mathbb{P}_k(E)\}.$$

The space is only implicitly defined and the actual shape of a $\boldsymbol{v} \in \boldsymbol{V}_k(E)$ in the interior of the element $E$ is not known a-priori. We can note that, in the case of a simplicial element $E$ with straight edges or faces, the space $\boldsymbol{V}_k(E)$ coincides with the classical Raviart-Thomas space $\mathbb{RT}_k(E)$ of order $k$ used in the mixed finite element method (Raviart and Thomas 1977). Furthermore, the local discrete space for the scalar trial $p$ and test $v$ fields of problem (1) is given by the classical approximation polynomial space

$$Q_k(E) := \{v \in Q(E) : v \in \mathbb{P}_k(E)\}.$$

The global spaces $\boldsymbol{V}_k(\Omega_h)$ and $Q_k(\Omega_h)$ are thus constructed by exploiting the continuity of the normal fluxes across each edge or face between two neighboring elements, i.e. $\boldsymbol{v} \cdot \boldsymbol{n}$ is single valued, and by the union of the local spaces, respectively. The space $Q_k(\Omega_h)$ is thus made of piecewise polynomials.

### Degrees of freedom

To make the objects defined on the spaces $\boldsymbol{V}_k(E)$ and $Q_k(E)$ computable, we need to introduce suitable degrees of freedom. These will be also useful in the actual computation of the bilinear forms in (1). For the space $\boldsymbol{V}_k(E)$, given an element $\boldsymbol{v}$ we introduce: *(i)* for each $e \in \mathcal{E}_h(E)$ the boundary moments as

$$\frac{1}{|h_e|} \int_e \boldsymbol{v} \cdot \boldsymbol{n} \, \widetilde{m}_d \, de \quad \forall \widetilde{m}_d \in \widetilde{\mathcal{M}_k}(e), d = 1, \ldots, \pi_k^{n-1},$$

and also *(ii)* the element moments of the divergence

$$\frac{h(E)}{|E|} \int_E \nabla \cdot \boldsymbol{v} \, m_d \, dE \quad \forall m_d \in \mathcal{M}_k(E) \setminus \mathcal{M}_0(E), d = 2, \ldots, \pi_k^n.$$

The space $\mathcal{M}_0(E)$ is clearly removed from the definition since gives null divergence, and it becomes useless in the definition of these degrees of freedom. Finally, we introduce the third set of degrees of freedom for $\boldsymbol{V}_k(E)$, namely *(iii)* the element moments

$$\frac{1}{|E|} \int_E \boldsymbol{v} \cdot \boldsymbol{m}_d^\oplus \, dE \quad \forall \boldsymbol{m}_d^\oplus \in \mathcal{G}_d^\oplus(E), d = 1, \ldots, \pi_{k-1}^n.$$

We note that, for application purposes, the *(i)* set of degrees of freedom for the velocity can be immediately used, without any post-processing, as an edge or face flux in the coupling with a transport or multi-phase flow problem.

Finally, the degrees of freedom for the space $Q_k(E)$ are just *(i)* the element moments

$$\frac{1}{|E|}\int_E v\, m_d\, dE \quad \forall m_d \in \mathcal{M}_k(E), d = 1, \dots, \pi_k.$$

A graphical representation of these degrees of freedom is given in Figure 2 for the special case of corner point cells.

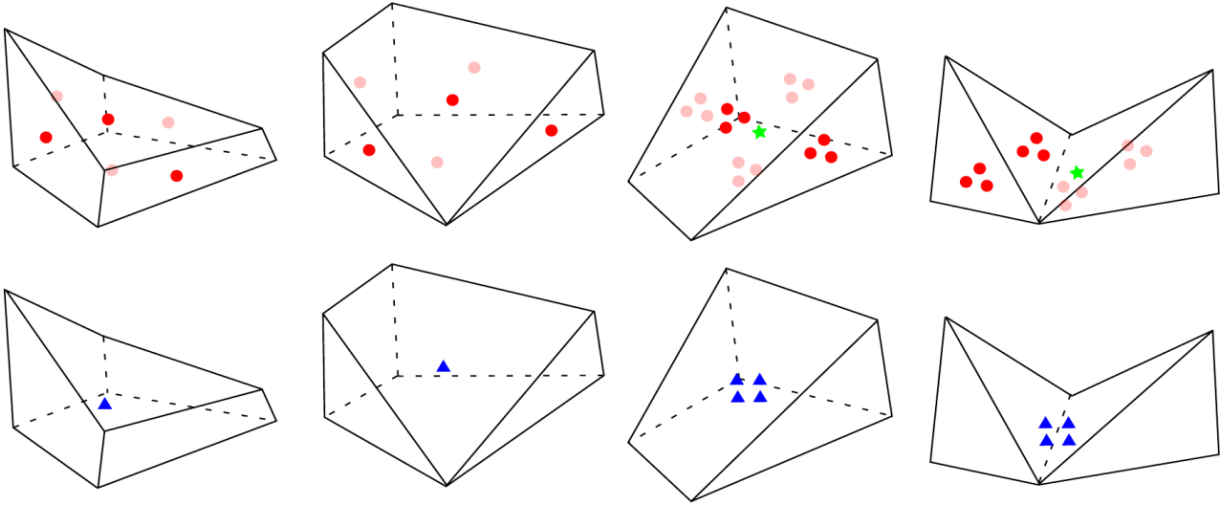

*Figure 2: On the top we indicate the dofs for the velocity and in the bottom for the pressure. The two cells on the left for $\boldsymbol{k} = \boldsymbol{0}$ and the two on the right for $\boldsymbol{k} = \boldsymbol{1}$. A red dot represents a (i) dof for $\boldsymbol{V}_k(\boldsymbol{E})$, a green star a (iii) dof for $\boldsymbol{V}_k(\boldsymbol{E})$, and a blue triangle a dof for $\boldsymbol{Q}_k(\boldsymbol{E})$.*

### Projection operators and approximation

As mentioned before, the actual shape of the basis functions for $\boldsymbol{V}_k(E)$ is not explicitly defined. However, the computation of the bilinear form $a$ requires their knowledge inside $E$. To overcome this issue, we introduce the projection operators $\Pi_0^d$ that maps $\boldsymbol{V}_k(E)$ to $[\mathbb{P}_d(E)]^n$

$$\int_E (\Pi_0^d \boldsymbol{v} - \boldsymbol{v}) \cdot \boldsymbol{m}\, dE = 0 \quad \forall \boldsymbol{m} \in [\mathcal{M}_d(E)]^n,$$

and $T_0^d \coloneqq I - \Pi_0^d$ as the orthogonal projection operator. It is possible to show that with the considered degrees of freedom of $\boldsymbol{V}_k(E)$ it is possible to compute the projection $\Pi_0^d$, for more details see (Dassi and Scacchi 2020).

Since the space $\widetilde{\mathbb{P}}_k(e)$ might not include polynomials up to order $k$ for $e \in \mathcal{E}_h(E)$, then also $\boldsymbol{V}_k(E)$ might not include polynomials up to order $k$ in $E$. It is possible to show that the constant polynomials are included in $\widetilde{\mathbb{P}}_k(e)$ and in $\boldsymbol{V}_k(E)$, too. However, to obtain optimal rate of convergence for the considered method we introduce the following extended space:

$$\boldsymbol{W}_k(E) \coloneqq \boldsymbol{V}_k(E) + [\mathbb{P}_k(A)]^n.$$

The local approximation of order $k$ of the bilinear form $a$, defined as $a_k : \boldsymbol{W}_k(E) \times \boldsymbol{W}_k(E) \to \mathbb{R}$, is thus given by

$$a(\boldsymbol{u}, \boldsymbol{v}) \approx a_k(\boldsymbol{u}, \boldsymbol{v}) \coloneqq a(\Pi_0^k \boldsymbol{u}, \Pi_0^k \boldsymbol{v}) + s(T_0^k \boldsymbol{u}, T_0^k \boldsymbol{v}).$$

The first part of $a_k$ is the consistency term, while the bilinear form $s$ is usually called stabilization. The former is computable in $\boldsymbol{V}_k(E)$ with the considered degrees of freedom. The latter is defined as the form $s : \boldsymbol{W}_k(E) \times \boldsymbol{W}_k(E) \to \mathbb{R}$ given by

$$s(\boldsymbol{u}, \boldsymbol{v}) \coloneqq \varsigma(E) \sum_{i,j=0}^{N_{dof}(E)} dof_i(\boldsymbol{u})\, dof_j(\boldsymbol{v}),$$

where $\varsigma$ is a parameter depending on the physical parameters $\mu K^{-1}$, on the element size and on the geometrical dimension $n$. Again, the bilinear form $s$ is computable in $V_k(E)$ with the introduced degrees of freedom. The other bilinear form $b$ involves only the divergence of the trial $u$. Thus, by considering the *(ii)* set of degrees of freedom of $V_k(E)$ we can immediately approximate it. The same applies for the boundary and source terms. The discrete problem we are solving is the following: *Find* $(q, p) \in V_k(\Omega_h) \times Q_k(\Omega_h)$ *such that*

$$\begin{cases} a_k(q,v) + b(v,p) = -(\overline{p}, v \cdot n)_{\partial_n \Omega_h} & \forall v \in V_k(\Omega_h) \\ b(q,v) = (f,v)_{\Omega_h} & \forall v \in Q_k(\Omega_h). \end{cases}$$

All the bilinear forms and functionals involve the computation of integrals, which are be done by suitable quadrature rules. We consider the quadrature rule introduced in \cite{}, that was already extended for the two-dimensional curved case (Beirao da Veiga, Russo and Vacca 2019) in two dimensions. Details on the application to the three-dimensional case will be shown in a forthcoming paper.

## Numerical Results

The present section is devoted to the discussion of some numerical tests, proposed to show and discuss the advantages and the applicability of the proposed technique. Two examples with a known analytical solution are included to validate the method and compare its performances with respect to standard approximation strategies, as well as two more applicative tests on more complex configurations, namely a domain cut by a listric fault and a corner point grid.

### Internal interface 2D

The first proposed example considers the two-dimensional square domain $\Omega = (0,1)^2$ shown in Figure 3 (left), where an inclusion $\Omega_2$ of radius $R = 0.45$ is present. The permeability of region $\Omega_1$ is $K_1 = I$, whereas the permeability in $\Omega_2$ is $K_2 = 0.01 I$. Viscosity is $\mu = 1$ in the whole domain. Continuity of the solution and flux balance are prescribed at the internal interface $\partial \Omega_2$, whereas boundary conditions on $\partial \Omega_1$ and forcing terms in $\Omega$ are set such that the exact pressure solution of equation (1) in $\Omega$ is $p_1(x,y) = K_2(x^2 + y^2) + R^2(I - K_2)$ in $\Omega_1$ and $p_2(x,y) = (x^2 + y^2)$ in $\Omega_2$, while $\mu q_i = -K_i \nabla p_i$, for $i = 1, 2$.

The problem is solved with both the here described method for curvilinear edges and with the standard MVEM.

The computational mesh, shown in Figure 3 (right), is generated with the following procedure: a quadrilateral mesh is built in $\Omega$ independently of the interface, which is subsequently superimposed, and the elements intersected by $\partial \Omega_2$ are cut into sub-elements of arbitrary shape not crossing it. This process is likely to generate badly shaped and elongated elements, as the ones in the region circled in the picture, which however can be easily handled by virtual-element based approaches [ (Berrone e and Borio 2017), (Mascotto 2018)], thus making the generation of the computational mesh an extremely simple task also for complex geometries with multiple intersecting interfaces.

For the MVEM for curvilinear edge elements, the geometry of the interface $\partial \Omega_2$ is kept unchanged, thus generating sub-elements with curvilinear edges, as can be seen in Figure 4 (right), whereas the interface is approximated by a linear interpolant connecting the intersection points between $\partial \Omega_2$ and mesh elements for the resolution with standard MVEM.

Let us introduce the following error indicators:

$$e_q = \sqrt{\sum_{E \in \Omega_h} |q - \Pi_0^k q_h|_E^2}, \quad e_p = \sqrt{\sum_{E \in \Omega_h} |p - p_h|_E^2},$$

where in the computation of $e_q$ we make use of the projected solution $\Pi_0^k q_h$ as proxy of the virtual solution $q_h$. We also denote the mesh-size parameter $h$, defined as $h = \frac{1}{N_E} \sum_{E \in \Omega_h} h(E)$. The convergence curves of the proposed error indicators with respect to the mesh-size parameter are reported in Figure 4, on the left for $e_p$ and on the right for $e_q$. In the picture, the curvilinear edge element MVEM is termed, for brevity, as *Curved* and the standard MVEM as *Straight*. Polynomial degrees ranging between 1 and 3 are considered and denoted as *deg k* in the figure legend. We can see that, as expected, only the *Curved* approach provides the optimal theoretical convergence trends for all the considered polynomial accuracy levels, whereas the *Straight* curves have maximum slope equal to two, independently of the polynomial degree, as the geometrical representation error of the curvilinear interface, proportional to $h^2$, dominates on the approximation error. This example also highlights the robustness of MVEM to badly shaped elements, curved elements and neighbouring cells with high aspect ratio.

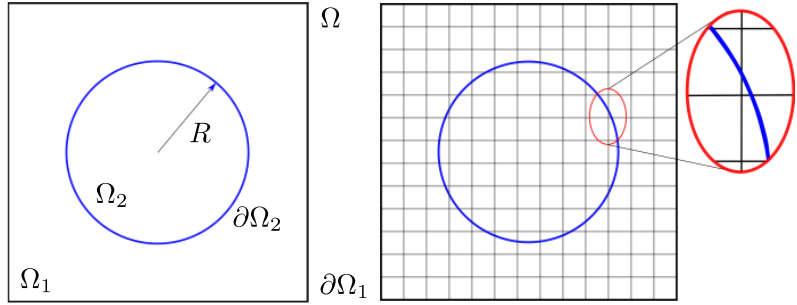

*Figure 3 Internal interface 2D. Domain description (left) and example computational mesh (right).*

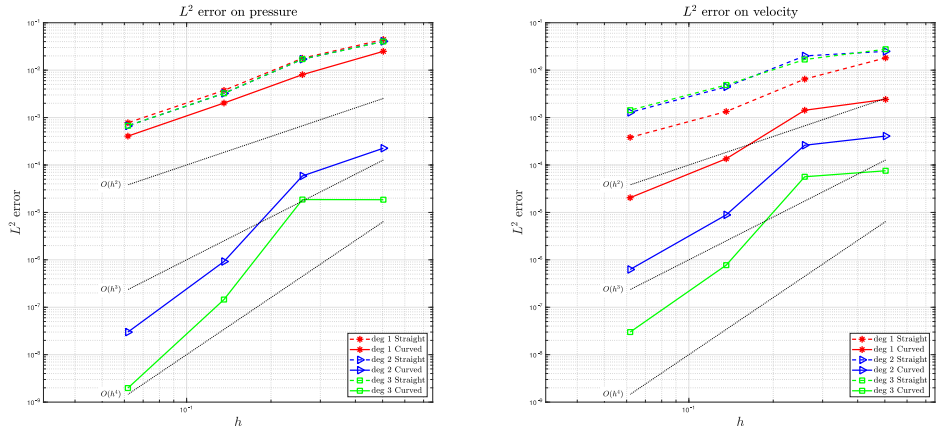

*Figure 4 Internal interface 2D. Convergence curves for pressure (left) and Darcy velocity (right). Method for curvilinear edges (Curved) compared to standard (Straight) approximation.*

## Curvilinear boundary 3D

The second example takes into account a 3D problem set on the domain represented in Figure 5: five of the six faces of the domain are planar and coincident with those of a unitary edge cube, whereas the top face Γ is curvilinear and defined by

$$\Gamma(x, y) = -\frac{1}{10} sin(\pi x) + 1.$$

Boundary conditions and forcing terms are prescribed such that the solution of problem (1) is given by

$$p(x, y, x) = \left(z + \frac{1}{10} sin(\pi x) - 1\right)^2,$$

with $\boldsymbol{K} = \boldsymbol{I}$ and $\mu = 1$. Null flux is imposed on the boundary and to recover the uniqueness of the solution we impose the exact pressure average. The geometry of the boundary Γ is exactly reproduced by the discrete mesh for the MVEM for curvilinear borders, whereas it is approximated by planar faces for standard MVEM. Convergence curves for the previous error indicators are shown, for this example, in Figure 6, where polynomial accuracy values ranging between 1 and 3 are considered and the *Curved* and *Straight* cases are compared. It can be seen that also for this example optimal convergence trends are achieved only by the curvilinear MVEM, being the error dominated by geometrical approximation in standard MVEM, thus bounding error decay to $O(h^2)$.

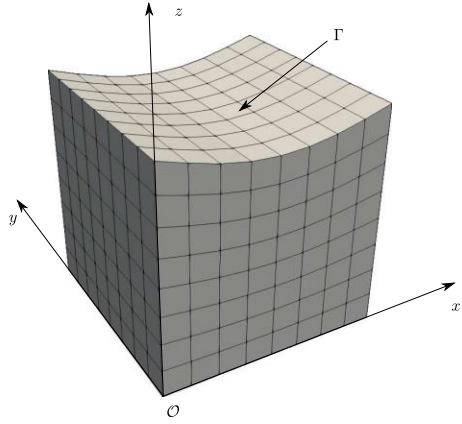

*Figure 5 Curvilinear boundary 3D. Domain and example of computational mesh.*

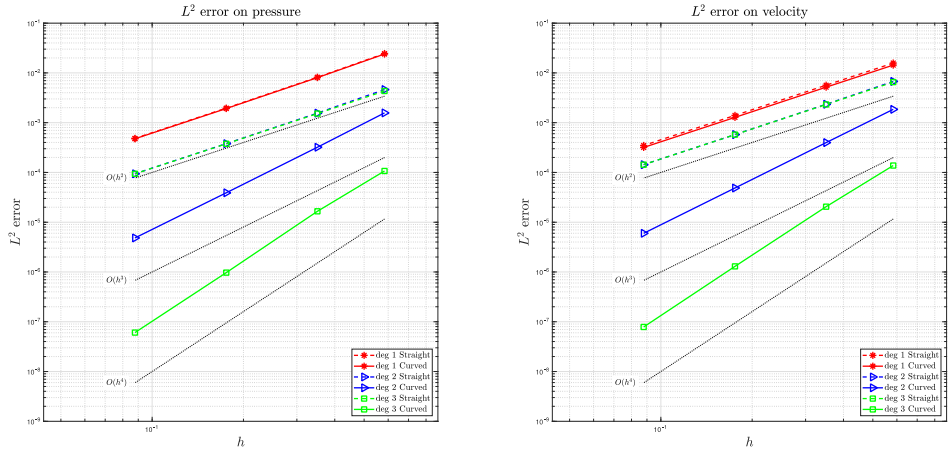

*Figure 6 Curvilinear boundary 3D. Convergence curves for pressure (left) and Darcy velocity (right). Method for curvilinear edges (Curved) compared to standard (Straight) approximation.*

### Listric fault

The last example takes into account a 2D rectangular domain $\Omega = [-1, 1] \times [-0.5, 0.5]$ with five interfaces $\Gamma_i$, $i = 1, \dots, 5$ that subdivide the domain into six subregions $\Omega_i$, $i = 1, \dots, 6$, as shown in Figure 7. The interfaces can be interpreted as a listric fault ($\Gamma_1$) and the broken horizons between sedimentary layers. The interfaces are defined as follows:

$\Gamma_1 = \{(x,y): y = -1.25\sqrt{(x+1.1)} + 1.01\}$,

$\Gamma_2 = \{(x,y): y = 0.25(x+1.1)^2 + 0.01\}$,

$\Gamma_3 = \{(x,y): y = 0.25(x+1.1)^2 - 0.21\}$,

$\Gamma_4 = \{(x,y): y = 0.5\sqrt{(x+1.1)} - 0.41\}$,

$\Gamma_5 = \{(x,y): y = 0.75\sqrt{(x+1.1)} - 0.01\}$.

Problem (1) is solved in $\Omega$ with the following parameters: $\mu = 1$, $\boldsymbol{K}_i = \xi_i \boldsymbol{I}$, with $\xi = (1, 0.01, 1, 1, 0.01, 1)$, and null forcing term. A null pressure is set on the top edge of the domain, a unitary pressure on the bottom edge, all other edges are instead insulated. An approximation order $k = 2$ is used. The obtained pressure distribution in $\Omega$ is presented in Figure 8, on a coarse mesh (left), as the one reported in Figure 7, and on a refined mesh (right): the behavior of the solution on the two meshes is in good agreement, and is also in agreement with the expected one. The steep change in pressure across the interfaces is clearly visible. We remark that, only for representation purposes, in Figure 8 interfaces are shown as piecewise linear and the solution is piecewise constant on mesh cells; the exact geometry of the interfaces and second order approximation for the solution are instead used in the computations.

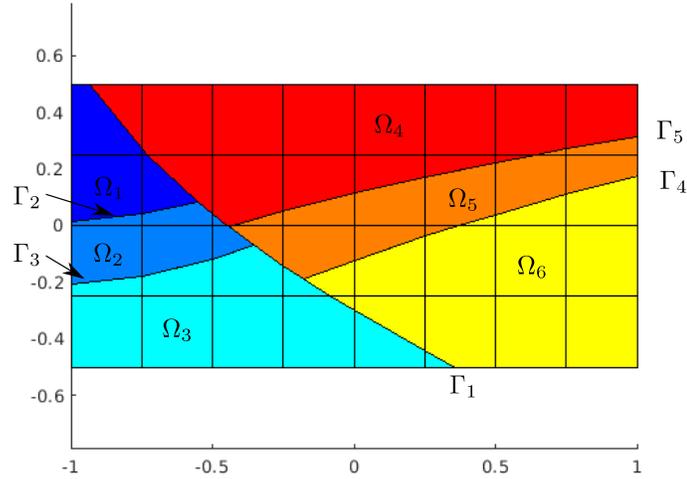

*Figure 7: Listric fault. Domain description*

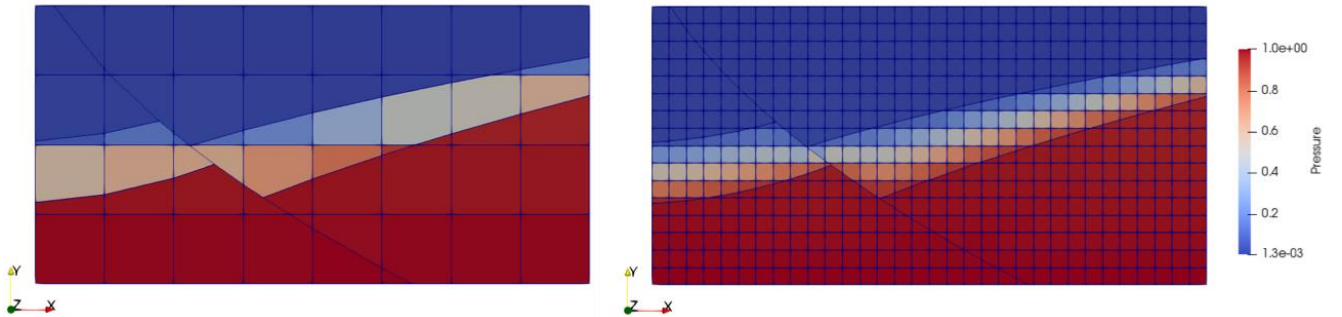

*Figure 8: Listric fault. Solution on coarse (left) and fine (right) meshes*

### Corner point mesh

The third proposed example takes into account a unit edge cube, composed of three layers of materials with different values of permeability. With reference to Figure 9, the lowest layer, coloured in blue, has permeability $\boldsymbol{K_1} = \boldsymbol{I}$, the middle layer, in gray, $\boldsymbol{K_2} = 0.01\boldsymbol{I}$, and the top layer, in red $\boldsymbol{K_3} = \boldsymbol{I}$. The two interfaces between the materials with different permeability, $\Gamma_1$ and $\Gamma_2$ are not planar, but bilinear, and the mesh is conforming to such interfaces. The mesh, indeed, is composed of cells with planar vertical faces and bi-linear top and bottom faces, with the exception of cell-faces lying on the external boundary of the domain. The equation of each bilinear surface, for $k = 1, ..., 9$ is defined by:

$$\sigma_k(u,v) = \boldsymbol{x}_A(1-u)(1-v) + \boldsymbol{x}_B u(1-v) + \boldsymbol{x}_C uv + \boldsymbol{x}_D(1-v)u$$

being $\boldsymbol{x}_A, \boldsymbol{x}_B, \boldsymbol{x}_C, \boldsymbol{x}_D$, with $\boldsymbol{x}_l \coloneqq (x_l, y_l, z_l)$, the coordinates of four points of each surface, placed at the intersections between the surface and the vertical edges of the domain, see Figure 9 (right), and:

$$u = \frac{x - x_A}{x_B - x_A}, \quad v = \frac{y - y_A}{y_D - y_A}.$$

A natural boundary condition $\bar{p} = 0$ is set on the top face of the domain, and a value $\bar{p} = 1$ is prescribed, instead, on the bottom face, all other faces having no-flux essential boundary conditions. Problem (1) is solved on such domain, with zero forcing term and $\mu = 1$, and approximation order $k = 2$. The obtained solution is reported in Figure 10, where the coloring is proportional to the computed pressure and two sections of the solution are also proposed. Figure 11 proposes a slice of the 3D solution on a plane normal to the $z$-axis and passing through the barycentre $(0.5, 0.5, 0.5)$ of the domain. The steep change of the solution across interfaces can be clearly seen. The trace on the plane of the interfaces is highlighted in red and is approximated by a piecewise linear interpolant only for graphical purposes, whereas the exact geometry is used for the computations. The method MVEM for curvilinear border elements can easily deal with corner point meshes, without introducing any geometrical representation error. This is a key property for underground flow simulations, where such kind of non-planar interfaces often arise in applications.

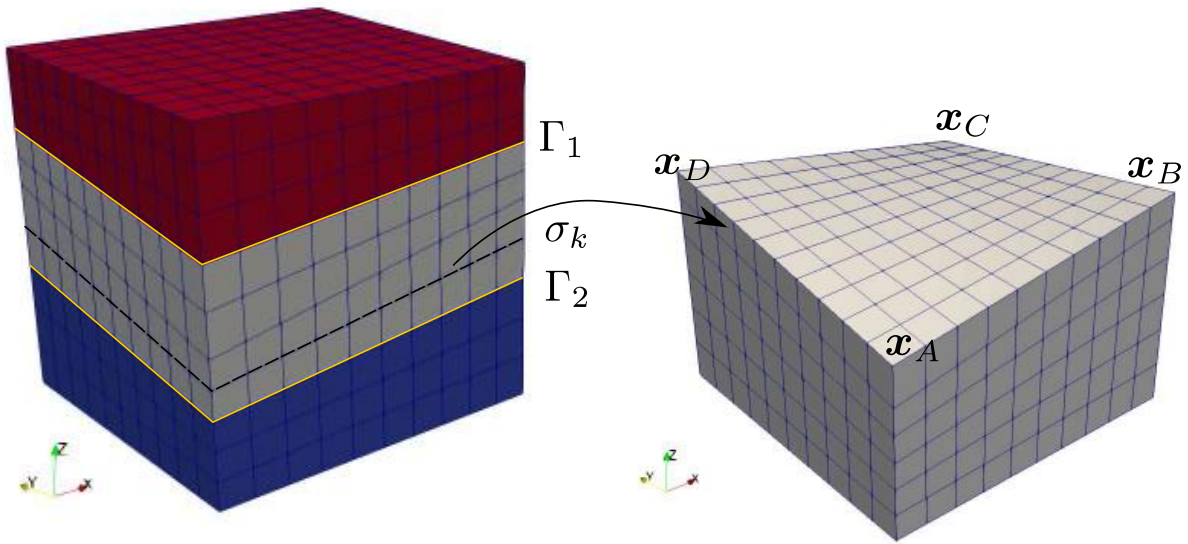

*Figure 9: Corner point mesh. Domain and material composition (left); mesh definition (right).*

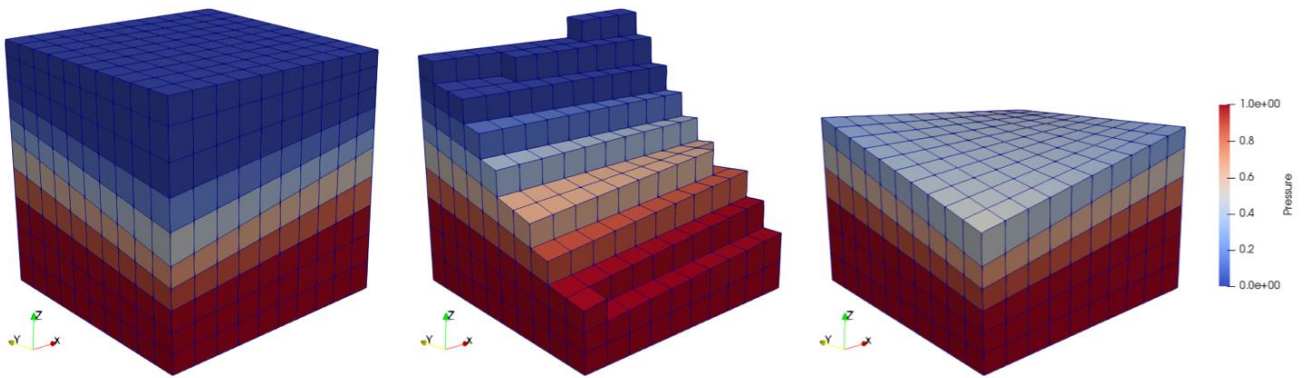

*Figure 10: Corner point mesh. Pressure solution (left) with two sliced views (middle and right)*

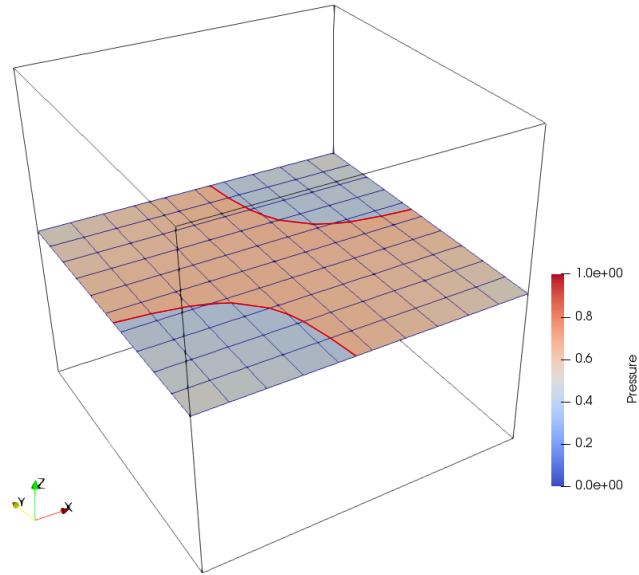

*Figure 11: Corner point mesh. Solution on a plane normal to the z-axis passing through domain barycentre*

## Conclusions

In this work we have introduced an innovative extension to the recent mixed-virtual element method for problem with curved interfaces. For high order approximation, a quality description of problem geometry can be of paramount importance in many real applications, where geometrical error dominates the one introduced by the numerical approximation scheme. The single-phase flow problem is kept in its mixed formulation so that both pressure and Darcy flux are directly computed avoiding, for the latter, a post-processing which might deteriorate its accuracy. We have introduced the finite dimensional spaces for both the pressure and velocity, being the latter virtual since the shapes of the basis functions are not prescribed a-priori. Following the virtual element method philosophy, we have introduced a projection operator and a stabilization bilinear form so that all the ingredients become computable. In the last part of the work, we have presented four test cases showing the validity of the proposed approach in both two and three-dimensions and considering geometries and challenges typical of underground applications. In all the cases we have obtained accurate results, making the proposed approach attractive for real life applications.

## Acknowledgements

The authors acknowledge financial support of INdAM-GNCS through project "Bend VEM 3d", 2020. Author S.S. also acknowledges the financial support of MIUR through project "Dipartimenti di Eccellenza 2018-2022" (Codice Unico di Progetto CUP E11G18000350001).